\documentclass{amsart}

\usepackage[normalem]{ulem}

\usepackage{lineno}

\usepackage{soul}

\usepackage{amssymb, amsmath, amsthm, hyperref, euscript, enumerate, color}
\usepackage[dvipsnames]{xcolor}
\usepackage{tikz-cd}

\usepackage{amscd}

\usepackage{bbm}
\usepackage{color,soul}
\usepackage[shortlabels]{enumitem}
\usepackage[english]{babel}

\DeclareMathOperator{\consum}{\#}

\theoremstyle{plain}

\newtheorem{theorem}[equation]{Theorem}

\newtheorem{corollary}[equation]{Corollary}

\newtheorem{proposition}[equation]{Proposition}

\newtheorem{lemma}[equation]{Lemma}

\newtheorem{remark}{Remark}

\newtheorem{thm}{Theorem}

\newtheorem{cor}[thm]{Corollary}

\newtheorem{exa}[thm]{Example}

\theoremstyle{definition}

\newtheorem{definition}[equation]{Definition}

\newtheorem{example}[equation]{Example}

\newcommand{\spcob}{\Spin(2, n - 1)_0\text{-pseudo-Riemannian cobordism}}

\newcommand{\R}{\mathbb{R}}

\newcommand{\Z}{\mathbb{Z}}

\newcommand{\Spin}{\textnormal{Spin}}

\newcommand{\So}{\textnormal{SO}}

\begin{document}

\author{Valentina Bais, Victor Gustavo May Custodio and Rafael Torres}

\title[Spin structures on pseudo-Riemannian cobordisms]{Existence results of $\Spin(2, n - 1)_0$-pseudo-Riemannian cobordisms.}

\address{Scuola Internazionale Superiori di Studi Avanzati (SISSA)\\ Via Bonomea 265\\34136\\Trieste\\Italy}

\email{\{vbais, vmaycust, rtorres\}@sissa.it}

%\subjclass[2020]{Primary 57K45, 57R55; Secondary 57R40, 57R52}

\maketitle

\emph{Abstract}: In this note, we study necessary and sufficient conditions for the existence of a Spin $(n + 1)$-dimensional cobordism that supports a non-singular and non-degenerate pseudo-Riemannian metric of signature $(2, n - 1)$, which restricts to a non-singular time-orientable Lorentzian metric on its boundary. The corresponding cobordism groups are computed. 

\tableofcontents

\section{Introduction and main results.}

%The main objects of study in this note are cobordisms, namely triples $(W;M_1,M_2)$ consisting of a smooth compact connected $(n+1)-$manifold $W$ with non-empty boundary $\partial W=M_1\sqcup M_2$. In particular, given two manifolds $M_1, M_2$,  we are interested in understanding what topological conditions should be imposed on them in order to admit a special kind of cobordism, where some geometric conditions are required to coexist. Such conditions are explained in the following definition of $\Spin(2, n - 1)_0$-pseudo-Riemannian cobordism. 

Reinhart \cite{[Reinhart]} and Sorkin \cite{[Sorkin]} determined necessary and sufficient topological conditions for a compact $(n + 1)$-manifold $W$ to admit a non-singular time-orientable Lorentzian metric $g^L$ inducing a Riemannian metric on its boundary $\partial W = M_1\sqcup M_2$. In particular, a pair $((W; M_1, M_2), g^L)$ satisfying such properties is known as an $\So(1, n)_0$-Lorentzian cobordism (Definition \ref{Definition Lorentzian Cobordism}) and the existence of the Lorentzian metric $g^L$ on $W$ has been completely characterized by Reinhart and Sorkin by means of the Euler characteristics of $W$ and of the closed $n$-manifolds $M_1$ and $M_2$. These objects have been studied extensively throughout the years by several mathematicians and mathematical physicists, including Chamblin \cite{[Chamblin]}, Geroch \cite{[Geroch]} Gibbons-Hawking \cite{[GibbonsHawking1]}, Reinhart \cite{[Reinhart]}, Sorkin \cite{[Sorkin]} and  Yodzis \cite{[Yodzis]}. 

An interesting question is whether two given topological and geometric properties can co-exist on a given cobordism $(W; M_1, M_2)$. Let us think for a second about the third-dimensional case to fix ideas. Any two closed oriented 3-manifolds are the boundary of some compact Spin 4-manifold since the third Spin-cobordism group is trivial \cite{[Milnor]}. The results of Reinhart and Sorkin that were just mentioned also guarantee the existence of a $\So(1, n)_0$-Lorentzian cobordism between any two closed oriented 3-manifolds. The latter cobordism, however, need not support a Spin-structure. Gibbons-Hawking \cite{[GibbonsHawking1]} showed that the Kervaire semi-characteristic (see Definition \ref{Definition Kervaire Semi}) of $\partial W = M_1\sqcup M_2$ is the only topological obstruction for a $\So(1, n)_0$-Lorentzian cobordism to admit a compatible structure of a Spin-cobordism. Smirnov-Torres \cite{[SmirnovTorres]} generalized their results to arbitrary dimensions and computed the corresponding $\Spin(1, n)_0$-Lorentzian cobordism groups. 

In this paper, we extend these results further into the pseudo-Riemannian realm and occupy ourselves with the study of the following objects. 

\begin{definition}\label{Definition Kleinian Cobordism}An $\So(2, n - 1)_0$-pseudo-Riemannian cobordism between closed smooth oriented $n$-manifolds $M_1$ and $M_2$ is a pair\begin{align}\label{SO Cobordism}((W; M_1, M_2), g)\end{align} that consists of

(A) a cobordism $(W; M_1, M_2)$,

(B.1) a non-singular indefinite metric $(W, g)$ of signature $(2, n - 1)$ such that 

(C.1) its restriction to the boundary $\partial W = M_1 \sqcup M_2$ gives rise to non-singular time-orientable Lorentzian metrics $(M_1, g^L_{M_1})$ and $(M_2, g^L_{M_2})$; please see Section \ref{Section Structure Groups} for an explanation on our notation. 

When the cobordism $(W; M_1, M_2)$ is a $\Spin$-cobordism, we call (\ref{SO Cobordism}) a $\Spin(2, n - 1)_0$-pseudo-Riemannian cobordism between $M_1$ and $M_2$. We say that $M_1$ is $\Spin(2,n-1)_0$-cobordant to $M_2$.
%define cobordisms and spin cobordisms?
\end{definition}

A canonical example of a $\Spin(2, n - 2)_0$-pseudo-Riemannian cobordism is the following. 

\begin{exa}
\label{Example A}
Let $(X, g)$ be a closed Riemannian $(n-2)$-manifold of dimension at least two, take the 2-disk $D^2$ with polar coordinates $(r, \theta)\in D^2$, and consider the compact pseudo-Riemannian $n$-manifold
\begin{align}\label{Cobordism Example}
(D^2, -dr^2 - r^2d\theta^2)\times (X, g)
\end{align}
with the indefinite product metric of signature $(2, n-2)$. Furthermore, if $X$ admits a $\Spin$-structure, then (\ref{Cobordism Example}) is a $\Spin(2, n - 2)_0$-pseudo-Riemannian cobordism with boundary the Lorentzian $(n - 1)$-manifold
\begin{align}
(S^1, - d\theta^2)\times (X, g).
\end{align}

\end{exa}

%The topic of Lorentzian and $\Spin(1, n)_0$-Lorentzian cobordisms, where the cobordism $W$ in (\ref{SO Cobordism}) admits a Lorentizan metric that restricts to a Riemannian metric on the boundary $\partial W = M_1\sqcup M_2$, has been studied extensively by several authors; see, for example, Chamblin \cite{[Chamblin]}, Gibbons-Hawking \cite{[GibbonsHawking1]}, Reinhart \cite{[Reinhart]}, Smirnov-Torres \cite{[SmirnovTorres]} and the references therein. 

While some partial results on $\Spin(2, 2)_0$-pseudo-Riemannian cobordisms have been obtained by Alty-Chamblin in \cite{[AltyChamblin]}, there is no systematic study of these objects available in the literature. The goal of this paper is to fill such a gap.

Let $M_1$ and $M_2$ be two closed smooth $n$-manifolds that are $\Spin$-cobordant. In our first main result, we provide an almost complete topological characterization for the existence of a $\spcob$ $((W;M_1,M_2),g)$; see Remark \ref{Remark Complete}. The Euler characteristic of a manifold $X$ is denoted by $\chi(X)$ and the Kervaire semi-characteristic of an odd-dimensional manifold $Y$ is denoted by $\hat{\chi}_{\Z/2}(Y)$; see Definition \ref{Definition Kervaire Semi}.

\begin{thm}
\label{Theorem Extension3} 
Let $\{M_1,M_2\}$ be closed smooth $\Spin$-cobordant $n$-manifolds.

$\bullet$ Suppose $n\not\equiv 1,5,7 \mod 8$. There exists a $\Spin(2, n - 1)_0$-pseudo-Riemannian cobordism $((W;M_1,M_2),g)$ if and only if
\begin{enumerate}
    \item $\chi(M_1)=\chi(M_2)=0$ for $n$ even
    \item $\hat\chi_{\Z_2}(M_1)=\hat\chi_{\Z_2}(M_2)$ for $n$ odd.
\end{enumerate}

$\bullet$ Suppose $n\equiv 1\mod 4$. If $\hat\chi_{\Z_2}(M_1)=\hat\chi_{\Z_2}(M_2)$, there exists a $\spcob$.

$\bullet$ Suppose $n\equiv 7 \mod 8$. There is a $\spcob$ without any further assumptions.

\end{thm}

The main ingredients of our proof of Theorem \ref{Theorem Extension3} are results of Atiyah \cite{[Atiyah]}, Frank \cite{[Frank]}, Hirzebruch-Hopf \cite{[HirzebruchHopf]}, Matsushita \cite{[Matsushita2]}, Thomas \cite{[Thomas1], [Thomas2], [Thomas6], [Thomas7], [Thomas3], [Thomas4], [Thomas5]} on the existence of 2-distributions of the tangent bundle of an oriented smooth manifold along with work of Gibbons-Hawking \cite{[GibbonsHawking1]}, Kervaire \cite{[Kervaire]}, Kervaire-Milnor \cite{[KervaireMilnor]}, Lusztig-Milnor-Peterson \cite{[LusztigMilnorPeterson]} and Smirnov-Torres \cite{[SmirnovTorres]} relating the Euler characteristic of a $\Spin$-cobordism with the Kervaire characteristic of its boundary. 

In Section \ref{Section Groups}, we define the corresponding  $\Spin(2, n - 1)_0$-cobordism groups, which we denote by $\Omega^{\Spin}_{2, n - 1}$. We build on  Milnor's computations of $\Spin$-groups in low dimensions \cite{[Milnor]} in order to determine them. Along with Theorem \ref{Theorem Extension3}, the task yields the following depiction of the cobordisms of Definition \ref{Definition Kleinian Cobordism} in terms of simple topological invariants. We also compute the cobordism groups.

\begin{thm}\label{Theorem Extension2}Let $\{M_1, M_2\}$ be closed $\Spin$ $n$-manifolds.

$\bullet$ If $n = 3$, there is a $\Spin(2, 2)_0$-pseudo-Riemannian cobordism $((W; M_1, M_2))$ if and only if $\hat{\chi}_{\Z/2}(M_1) = \hat{\chi}_{\Z/2}(M_2)$. There is a group isomorphism\begin{align}\label{3D Isomorphism}\Omega^{\Spin_0}_{2, 2} \rightarrow \Z/2.\end{align}

$\bullet$ If $n = 4$, there is a $\Spin(2, 3)_0$-pseudo-Riemannian cobordism $((W; M_1, M_2), g)$ if and only if $M_1$ and $M_2$ have trivial Euler characteristic and the same signature, i.e. $\chi(M_1)=0=\chi(M_2)$ and $\sigma(M_1)=\sigma(M_2)$. There is a group isomorphism\begin{align}\label{4D Isomorphism}\Omega^{\Spin_0}_{2, 3} \rightarrow \Z.\end{align}

$\bullet$  If $n = 6$, there is a $\Spin(2, 5)_0$-pseudo-Riemannian cobordism $((W; M_1, M_2), g)$ if and only if $\chi(M_1) = 0 = \chi(M_2)$. There is a group isomorphism\begin{align}\Omega^{\Spin_0}_{2, 5}\rightarrow \{0\}.\end{align}

$\bullet$  If $n = 7$, there is a $\Spin(2, 6)_0$-pseudo-Riemannian cobordism $((W; M_1, M_2), g)$ without any further assumptions. There is a group isomorphism\begin{align}\Omega^{\Spin_0}_{2, 6}\rightarrow \{0\}.\end{align}

\end{thm}

The last result to be presented in this introduction contains a myriad of  examples of the cobordisms of Definition \ref{Definition Kleinian Cobordism}; cf. \cite[Corollary G]{[SmirnovTorres]}. 

\begin{thm}\label{Theorem Extension4}
Let $n\geq 5$; let $\hat W$ be a compact smooth $(n+1)$-manifold that admits a $\text{Spin}(n+1)$-structure, and let $M=\partial W$ be its non-empty boundary. For any finitely presented group $G$, there exists a closed smooth $n$-manifold $M(G)$ such that $\pi_1 M(G)=G$ and a $\text{Spin}(2,n-1)_0$- pseudo-Riemannian cobordism $((W;M,M(G)),g)$. 
\end{thm}

The proof of Theorem \ref{Theorem Extension4} follows from Theorem \ref{Theorem Extension3} and a well-known argument to construct closed high-dimensional stably-parallelizable $n$-manifolds with prescribed fundamental group \cite[Theorem A]{[JohnsonWalton]}. %A proof of Theorem \ref{Theorem Extension4} is obtained via a slight modification to an argument of Smirnov-Torres \cite[Proof of Corollary G]{[SmirnovTorres]}, and it is therefore omitted.

We have organized the paper as follows. In Section \ref{background} we present some background material for the convenience of the reader. It includes a description of the topological constructions that are used in the paper as well as background existence results on indefinite metrics. A discussion on the co-existence of Spin structures and indefinite metrics on our cobordisms can be found in Section \ref{Section Structure Groups}. The cobordism groups are defined in Section \ref{Section Groups}. Section \ref{Section Comparison} contains a comparison between the cobordisms of Definition \ref{Definition Kleinian Cobordism} and Lorentzian cobordisms. The proofs of our main results are given in Section \ref{Section Proofs}. For background results, the reader is directed to Atiyah \cite{[Atiyah]}, Chamblin \cite{[Chamblin]}, Gibbons-Hawking \cite{[GibbonsHawking1]}, Milnor \cite{[Milnor]}, O'Neill \cite{[ONeill]}, Reinhart \cite{[Reinhart]}, Smirnov-Torres \cite{[SmirnovTorres]}, Steenrod \cite{[Steenrod]}, Stong \cite{[Stong]}, Thom \cite{[Thom]}, Thomas \cite{[Thomas1], [Thomas2], [Thomas3], [Thomas4], [Thomas5]}. 

All manifolds in this paper are assumed to be $C^{\infty}$-smooth and Hausdorff. All pseudo-Riemannian metrics in this paper are assumed to be time-orientable and non-degenerate.

\section{Background results}\label{background}

\subsection{Kervaire semi-characteristic and Spin-structures} The following fundamental invariant of odd-dimensional manifolds was introduced by Kervaire \cite{[Kervaire]}.

\begin{definition}\label{Definition Kervaire Semi}
The Kervaire semi-characteristic of a closed $(2k + 1)$-manifold $M$ is
\begin{align}
\hat{\chi}_{\Z/2}(M) = \sum_{i = 0}^k b_i(M; \Z/2) \mod 2
\end{align}
where $b_i(M;\Z/2)$ denotes the $i^{th}$ Betti number of $M$ with $\Z/2$-coefficients. %More generally, given a field $F$, we define
%\begin{align}
 %   \hat{\chi}_{F}(M) \equiv \sum_{i = 0}^k \text{dim} \,H_i(M;F) \mod 2.
%\end{align}
\end{definition}

There is a relation between the Euler characteristic of an even-dimensional manifold with boundary and the Kervaire semi-characteristic of its boundary in the presence of a $\Spin$-structure as observed by Geiges \cite{[Geiges]}, Gibbons-Hawking \cite{[GibbonsHawking1]}, Kervaire \cite{[Kervaire]}, Kervaire-Milnor \cite{[KervaireMilnor]}, Lusztig-Milnor-Peterson \cite{[LusztigMilnorPeterson]} and Smirnov-Torres \cite{[SmirnovTorres]}. 

\begin{theorem}\label{Theorem Relation}Let $W$ be a compact even-dimensional manifold with non-empty boundary $\partial W \neq \emptyset$. The identity \begin{align}\chi(W) + \hat{\chi}_{\Z/2}(\partial W) \equiv 0 \mod 2\end{align} holds provided either

\begin{itemize}
\item $W$ is stably-parallelizable \cite{[Kervaire]} or
\item $W$ is $2q$-dimensional for $q\not\equiv 0 \mod 4$ and admits a $\Spin$-structure \cite[Lemma 8.2.13]{[Geiges]}, \cite[Lemma 6]{[SmirnovTorres]}.
\end{itemize}
\end{theorem}

The value of the Kervaire semi-characteristic is independent of the choice of field of coefficients for manifolds that admit a $\Spin$-structure \cite{[LusztigMilnorPeterson]}.

\subsection{Double of a compact manifold and its Euler characteristic and Kervaire semi-characteristic}

Let $W$ be a compact oriented $n$-manifold with non-empty boundary $\partial W = M$. The double $2W$ of $W$ is the closed smooth oriented $n$-manifold\begin{align}\label{Double}2W = W \cup_M \overline{W},\end{align} where corners have been smoothed out. It can also be described as the boundary\begin{align}\label{Double Boundary}2W = \partial (W\times [0, 1]).\end{align}The computation of the Euler characteristic of (\ref{Double}) and its signature, whenever it is defined, are immediate and we record them. 

\begin{lemma}\label{Lemma Characteristic Double}The Euler characteristic of the double (\ref{Double}) is\begin{align}\label{Characteristic Double}\chi(2W) = 2\chi(W) - \chi(\partial W).\end{align} 

%Every Stiefel-Whitney and Pontrjagin number of $2W$ vanishes.

Suppose $W$ is $4q$-dimensional. The signature of $2W$ satisfies\begin{align}\sigma(2W) = 0.\end{align}

\end{lemma}

%The vanishing of the Stiefel-Whitney and Pontryagin numbers follows from the double of $W$ being a boundary (\ref{Double Boundary}); see \cite{[Kosinski], [Thom]}. 

A result of Zadeh \cite{[MEZadeh]} is useful for the computation of the Kervaire semi-characteristic of (\ref{Double Boundary}).

\begin{proposition}\label{Proposition Vanishing SCharacteristic}

    Let $W$ be a $\Spin$ compact $2q-1$ oriented manifold with $q\not \equiv 0\mod 4$. The Kervaire semi-characteristic of the double $2W$ satisfies 
    \begin{align}
        \hat{\chi}_{\Z/2}(2W) \equiv 0 \mod 2.
    \end{align}
\end{proposition}

The proof of Proposition \ref{Proposition Vanishing SCharacteristic} follows from Theorem \ref{Theorem Relation} and the fact that $2W$ is the boundary of the Spin-manifold $W\times [0,1]$.

\subsection{Existence of 2-distributions}\label{Section 2-distributions} Let us now discuss existence results of subbundles of the tangent bundle of a smooth manifold that are directly related to the existence of pseudo-Riemannian metrics.

\begin{definition}Let $W$ be a smooth oriented $n$-manifold. A 2-distribution $V$ is a non-singular field of tangent 2-planes, i.e. an oriented 2-plane sub-bundle of $TW$.
\end{definition}

\begin{remark} An oriented rank 2 vector bundle over a manifold is determined up to isomorphism by its Euler class {\cite{[Steenrod]}. It follows that a 2-distribution on $W$ with trivial Euler class satisfies $\text{Span}\;W\geq 2$}, where $\text{Span }W$ is the maximum number of everywhere linearly independent and nowhere vanishing vector fields on $W$.
\end{remark}

The following result collects several foundational theorems on the existence of 2-distributions on oriented and $\Spin$-manifolds due to Atiyah \cite{[Atiyah]}, Frank \cite{[Frank]}, Hirzebruch-Hopf \cite{[HirzebruchHopf]}, Matsushita \cite{[Matsushita2]} and Thomas \cite{[Thomas1], [Thomas2], [Thomas3], [Thomas4], [Thomas5], [Thomas6], [Thomas7]}. 

\begin{theorem}\label{Theorem Main} Let $X$ be a closed smooth oriented $n$-manifold.

$\bullet$ \cite{[Frank]}, \cite[4.5]{[HirzebruchHopf]}, \cite[Theorem 2]{[Matsushita1]}, \cite[Theorem 3.1]{[Atiyah]}. Suppose $n \equiv 0 \mod 4$ and that the signature satisfies $\sigma(X) = 0$. There is a 2-distribution on $X$ if and only if\begin{center}$\chi(X) \equiv 0 \mod 4$.\end{center}

$\bullet$ \cite[Theorem 1.3]{[Thomas3]}. Suppose $n \equiv 1 \mod 4$ and that $X$ admits a $\Spin$-structure. There is a 2-distribution on $X$ with Euler class $2v$ for each class $v\in H^2(X; \Z)$ if and only if
\begin{align*}
w_{n - 1}(X) = 0 \;\text{ and }\;\hat{\chi}_{\Z/2}(X) \equiv 0 \mod 2. 
\end{align*}
Moreover, if $X$ admits a 2-distribution, then $\hat\chi_\R(X)\equiv0 \mod 2$ \cite[Theorem 4.1]{[Atiyah]}.

$\bullet$ \cite[Theorem 1.1, Corollary 1.2]{[Thomas3]}. If $n \equiv 3 \mod 4$, then there is a 2-distribution on $X$ with Euler class $2v$ for each class $v\in H^2(X; \Z)$. 

\end{theorem}

Theorem \ref{Theorem Main} is a key ingredient in the proofs of our main results. The following technical lemma is used to build a $\Spin(2, n-1)_0$-pseudo-Riemannian cobordism whenever a given $\Spin$-cobordism $(W; M_1, M_2)$ satisfies $\chi(W)=\chi(M_1)=\chi(M_2)=0$ and $\text{Span}(W) \geq 2$.

\begin{lemma}\label{vec}
  Let $(W; M_1, M_2)$ be a cobordism between smooth closed orientable $n$-manifolds such that all of these manifolds have trivial Euler characteristic: $$\chi(W)=\chi(M_1)=\chi(M_2)=0.$$ Suppose furthermore that there exist two everywhere linearly independent vector fields $X, Y \in \mathfrak{X}(W)$. Then we can find two everywhere linearly independent vector fields $\tilde X, \tilde Y \in \mathfrak{X}(W)$ such that $\tilde X |_{\partial W}$ is the exterior normal to the boundary $\partial W$.
\end{lemma}
\begin{proof}
Fix a Riemannian metric $g$ on $W$. Let $X, Y \in \mathfrak{X}(W)$ be as above and $\nu_i$ a tubular neighbourhood of $N_i$ inside $W$ for $i=1,2$ respectively. In the following, we will make use of the identification $ N_i \times [0,1] \cong \nu_i$, with $N_i \times \{1\}$ corresponding to the submanifold $N_i \subset \nu_i$. Let $\bar X \in \mathfrak{X}(W)$ be a vector field satisfying the following conditions:
\begin{enumerate}
\item $\bar X |_{W \setminus (\nu_1 \sqcup \nu_2)}=X|_{W \setminus (\nu_1 \sqcup \nu_2)}$, i.e. $\bar X$ coincides with $X$ outside the disjoint union of the collar neighbourhoods of the boundary components of $W$;
\item $\bar X|_{\partial W}$ is the outward-pointing normal vector field;
\item $\bar X$ has finitely many singular points in $\nu_1 \sqcup \nu_2$.
\end{enumerate}
Notice that such a vector field can always be built out of $X$.
Up to composing with a self-diffeomorphism of $W$, we can assume that all the singular points of $\bar X$ are contained in the disjoint union of two small balls $B_i \subset \nu_i$ for $i=1,2$ (see \cite{[Milnor2]}, Chapter 4). 

Our first aim is building out of $Y$ a new vector field $\bar Y \in \mathfrak{X}(W)$ with the following properties:
\begin{enumerate}
    \item $\bar Y |_{W \setminus (\nu_1 \sqcup \nu_2)}=Y|_{W \setminus (\nu_1 \sqcup \nu_2)}$, i.e. $\bar Y$ coincides with $Y$ outside the disjoint union of the collar neighbourhoods of the boundary components of $W$;
    \item $\bar X, \bar Y$ are everywhere linearly independent in $W \setminus \text{Int}( B_1 \sqcup B_2) \cong W \setminus \{ \text{2 points}\}$.
\end{enumerate}
Let's build $\bar Y$ in the tubular neighbourhood $\nu_1$, the procedure in $\nu_2$ will be analogous. Without loss of generality, we can suppose $\bar X$ to be unitary outside $W \setminus \text{Int}( B_1 \sqcup B_2) \cong W \setminus \{ p_1, p_2\}$ with respect to the metric $g$, where $p_i =(x_i,t_i)\in \nu_i$ for $i=1,2$ respectively. The vector field $\bar X$ defines an isotopy between the linear span of $X|_{N_1 \times \{0\}}$ inside $TW|_{N_1 \times \{0\}}$ and the one of $\bar X|_{N_1 \times \{t_i\}}$ inside $TW|_{N_1 \times \{
t_i \}}\cong TW|_{N_1 \times \{0\}}$. Such isotopy can be extended to an ambient isotopy of the total space of the bundle $TW|_{N_1 \times \{0\}}$ consisting of vector bundle isomorphisms (by adapting the proof of Hirsch's Isotopy extension theorem in  \cite{[Hirsch]}). We can use such isotopy to define $\bar Y$ on $N_1 \times [0, t_i)$. Since $\chi(N_1)=0$, Poincaré- Hopf's theorem (see \cite{[Milnor2]}, Chapter 6) allows us to define $\bar Y$ globally also on $N_1 \times \{t_i\}$ and conclude by extending it to all $\nu_1$ as done before.

Once we have built such $\bar Y$, we can suppose (up to composing with a self-diffeomorphism of $W$) that $\bar X, \bar Y$ are linearly independent outside a little ball $B \subset \nu_1$. Moreover, we can assume without loss of generality that $\bar X, \bar Y$ are orthonormal with respect to the fixed Riemannian metric. Poincaré-Hopf's theorem together with the vanishing of the Euler characteristic of $W$ implies that the map
\begin{align*}
    \partial B \to S^n
\end{align*}
sending each point $x \in \partial B$ to $\bar X(x)$ is of zero degree and is therefore null-homotopic. Hence, we can suppose without loss of generality that the restriction $\bar X|_{\partial B}$ is the constant vector field $\frac{\partial}{\partial t}$ under the identification $\nu_1 \cong N_1 \times [0,1]$. In this way, $\bar X$ can be constantly extended in the interior of $B$, defining a global vector field $\tilde X \in \mathfrak{X}(W)$. 

Moreover, $\bar Y|_{N_1 \times \{t\}}$ can be seen as a vector field on $N_1 \setminus \{point\}$ whenever $(N_1 \times \{t\}) \cap B \neq \emptyset$. Since $\chi(N_1)=0$, Poincarè-Hopf's theorem allows us to suppose that $\bar Y |_{\partial B}$ is constantly equal to a fixed vector in the tangent space of $N_1$ and hence we are able to extend it to a global vector field $\tilde Y \in \mathfrak{X}(W)$ which is everywhere linearly independent with respect to $\tilde X$.
\end{proof}

\subsection{Existence of indefinite metrics}\label{Section Indefinite Metrics} We now recall some basic definitions and existence results of pseudo-Riemannian metrics with non-trivial signature on smooth manifolds. It is well known that the existence of a Lorentzian metric on a smooth manifold is equivalent to the existence of a nowhere vanishing vector field, i.e. a nowhere vanishing section of its tangent bundle \cite{[ONeill]}, \cite[Lemma 1]{[SmirnovTorres]}. Sub-bundles of the tangent bundle of a smooth manifold yield other indefinite metrics \cite{[ONeill], [Steenrod]}, and the role in this note of the existence results on 2-distributions that were described in Section \ref{Section 2-distributions} is explained in the following lemma; cf. \cite{[Steenrod]}.

\begin{lemma}\label{Lemma 2Distribution}Let $W$ be a smooth $n$-manifold. There is a pseudo-Riemannian metric $(W, g)$ of signature $(p, q)$ if and only if there is a decomposition of the tangent bundle\begin{align}\label{eqn: decomp pseudo-riem metric}TW = \xi \oplus \eta,\end{align}where $\xi$ and $\eta$ are vector sub-bundles of rank $p$ and $q$, respectively. Moreover, these sub-bundles can be chosen such that $\xi$ is time-like and $\eta$ is space-like.

In particular, there is a non-singular indefinite metric of signature $(2, n - 2)$ on a smooth $n$-manifold $W$ if and only if there is a 2-distribution $V\subset TW$.
\end{lemma}

The vector sub-bundle $\xi$ from the decomposition (\ref{eqn: decomp pseudo-riem metric}) is called a time-like sub-bundle of maximal rank. A pseudo-Riemannian manifold $(W,g)$ of signature $(p,q)$ is called time-orientable if there exists an orientable time-like sub-bundle of maximal rank. Here and throughout, all pseudo-Riemannian metrics are assumed to be time-orientable.

The following result will be used to construct a non-singular pseudo-Riemannian metric of signature $(2, n - 2)$ on the double $2W$ whenever there is such a metric on $W$, which restricts to a Lorentzian metric on $\partial W$.

\begin{theorem}\label{Theorem-gluemetric} Let $W$ be a smooth $n$-manifold with non-empty boundary $\partial W$. The following statements are equivalent.\begin{enumerate}
\item There is an indefinite metric $(W, g)$ of signature $(2, n - 2)$ such that its restriction to the boundary $g|_{\partial W}$ is a Lorentzian metric;
\item There is a rank 2 sub-bundle $\xi \subset TW$ and a line bundle $\eta'\subset \xi|_{\partial W}$ such that $\eta'$ is transversal to the zero section $\partial W$ of $\xi|_{\partial W}$;
\item There is a rank 2 sub-bundle $\xi \subset TW$ and a nowhere vanishing section $\nu:\partial W\rightarrow \xi|_{\partial W}$ that is everywhere outward pointing.

\end{enumerate}

\end{theorem}

\begin{proof}
Suppose the first item holds. Lemma \ref{Lemma 2Distribution} implies the existence of a decomposition
\begin{align*}
TW=\xi \oplus \eta 
\end{align*}
where $\xi, \eta \subset TW$ are sub-bundles of rank $2$ and $n-2$ respectively and such that $g$ is negative-definite on $\xi$ and positive-definite on $\eta$. Since $g':= g\bigr \vert_{\partial W}$
 is Lorentzian by hypothesis, Lemma \ref{Lemma 2Distribution} implies the existence of a further decomposition 
 \begin{align*}
 T\bigr(\partial W\bigr)=\tilde \xi \oplus \tilde \eta
 \end{align*}
 
 where $\tilde \xi$ is a time-like rank $1$ vector sub-bundle of $T\bigr(\partial W\bigr)$, while $\tilde \eta$ can be chosen to be space-like.  
 
 We will prove that the map 
 \begin{align*}
 p':\tilde \xi \to \xi\bigr\vert_{\partial W}
 \end{align*}
 given by the composing the inclusion of the bundle $\tilde \xi$ in $TW$ with the projection 
 \begin{align*}
 p:TW=\xi \oplus \eta \to \xi.
 \end{align*}
 is a bundle monomorphism.  
 
 Let $x\in \partial W$ and consider a non-zero vector $v\in \tilde\xi_x$. We will show that $w:=p'(v)\in \xi_x$ is then also non-zero.  Indeed, suppose that $w\in \xi_x$ is the zero vector. Then (after composing with the inclusion $\tilde \xi_x \subset TW$) we would have $v\in \eta_x$ and thus $g|_{T_xW}(v,v)>0$ - a contradiction. 
 
 Then, the image of $\tilde \xi$ under the bundle map $p': \tilde \xi \to \xi|_{\partial W}$ defines a line sub-bundle $\xi_1\subset \xi\bigr\vert_{\partial W}$, yielding a further decomposition 
 \begin{align*}
 \xi\bigr\vert_{\partial W}=\xi_1\oplus \xi_2.
 \end{align*}
 
 In order to show item 2, it will be enough to prove the transversality condition $\xi_2 \pitchfork \partial W$.
 
 Suppose by contradiction that there exists a point $x \in \partial W$ such that $(\xi_2)_x \subset T_x \partial W$. This implies that $(\xi_2)_x=<v>$ is generated by a nonzero vector $$v=v_1+v_2\in \tilde \xi_x \oplus \tilde \eta_x$$ where $v_1 \in \tilde \xi_x$ and $v_2 \in \tilde \eta_x$. Since $p'$ is an isomorphism and $v \notin (\xi_1)_x$, we have that $v_1=0$, while the fact of $v$ being time-like implies that also $v_2$ is trivial, a contradiction.

 Suppose now that item $2$ holds true. Let $\hat{n}:\partial W \to TW\bigr\vert_{\partial W}$ be an outward-pointing vector field. For all $x\in \partial W$ there is a decomposition 
 \begin{align*}
      T_xW=T_x\bigr(\partial W\bigr)\oplus \text{Span}\bigr(\hat n(x)\bigr).
 \end{align*}
and from the decomposition $T_x W= T_x(\partial W) \oplus \eta'_x$ we have $$\hat{n}(x)=v_x+w_x$$ where $v_x\in T_x\bigr(\partial W\bigr)$ and $w_x\in \eta'_x$. Note that the condition of $\hat n$ being outward pointing trivially implies that $w_x\neq 0$ for all $x\in \partial W$. In particular, the map 
$$n:x\mapsto w_x$$ defines a nowhere vanishing section of $\xi\bigr\vert_{\partial W}$. Since $w_x\notin T_x\bigr(\partial W\bigr)$ at every point by construction, we have that $n$ is either everywhere outward pointing or everywhere inward pointing. In the latter case we can simply consider $-n$ as the desired section and prove item 3.

Suppose finally that item 3 holds. We can use the given $2$-distribution $\xi$ to construct a pseudo Riemannian metric $g$ on $W$ of signature $(2, n-2)$ with $\xi \subset TW$ being time-like.  We consider also a rank 2 space-like sub-bundle $\eta\subset TW$ complementary to $\xi$. The everywhere outward-pointing vector field $n:\partial W\to \xi\bigr\vert_{\partial W}$ spans a trivial sub-bundle $\xi_1\subset \xi\bigr\vert_{\partial W}$. Hence, we have a decomposition $\xi\bigr\vert_{\partial W}=\xi_1\oplus \xi_2$. Clearly, we may assume that $\xi_2 \subset T\bigr(\partial W\bigr)$. We thus have a decomposition
\begin{align*}
    T\bigr(\partial W\bigr)=\xi_2 \oplus \eta\bigr\vert_{\partial W};
\end{align*}
with $\xi_2$ time-like and $\eta|_{\partial W}$ space-like. Such a decomposition allows us to get the desired Lorentzian metric on $\partial W$ and conclude. 
\end{proof}

At this point, we are able to formulate the following consequence of Lemma \ref{vec} and Theorem \ref{Theorem-gluemetric} that will be useful for our purposes. 

\begin{corollary}\label{cor}
Let $(W;M_1,M_2)$ be a $\Spin$-cobordism. The following statements are equivalent.

$\bullet$ There is an indefinite metric $g$ such that $((W;M_1,M_2),g)$ is a $\Spin(2,n-1)_0$-pseudo-Riemannian cobordism and $g$ is induced by a 2-distribution of trivial Euler class;

$\bullet$ There exist everywhere linearly independent vector fields $X,Y \in \mathfrak{X}(W)$ and $$\chi(W)=\chi(M_1)=\chi(M_2)=0.$$

\end{corollary}

We will make use of the necessary and sufficient conditions for a closed $n$-manifold $X$ to admit two everywhere linearly independent vector fields, with the purpose of building a $\spcob$ out of a $\Spin$-cobordism as in Corollary \ref{cor}. These conditions have been obtained by Thomas \cite[Table 2, pp. 652]{[Thomas5]} and are summarized in Table 1.

\begin{table}[h]
    \centering
    \begin{tabular}{c|c} \hline \hline
     \shortstack{$n=\dim X$ \\ $ \;$}  &  \shortstack{$\;$  \\Necessary and sufficient
        \\ conditions for $\text{Span}\; X\geq2$} \\ \hline 
        & \\[-1em]
        $n\equiv 0\mod 4$ & $\sigma(X)\equiv 0 \mod 4$, $\chi(X)=0$ \\ & \\[-1em]
        $n\equiv 1\mod 4$ & $w_{n-1}(X)=0$, $\hat\chi_{\mathbb R}(X)=0$ \\ & \\[-1em]
        $n\equiv 2\mod 4$ & $\chi(X)=0$ \\ & \\[-1em]
        $n\equiv 3\mod 4$ & $\text{Span}\;X\geq 2$ is always true \\ & \\ \hline 
    \end{tabular}
    \caption{ }
    
    \label{Thomas5_table}
\end{table}

\begin{proposition} \label{Proposition IDouble}
Let $W$ be a compact orientable $n$-manifold with nonempty boundary such that $\chi (\partial W)=0$.

$\bullet$ Suppose that $n\not \equiv 3 \mod 4$  and $\chi (W)=0$. If $n \equiv 1 \mod 4$, assume further that $W$ admits a $\Spin$-structure.  There are non-singular Lorentzian metrics $(2W, g^L)$ and $(W, g^L_W)$ such that the latter restricts to a Riemannian metric on $\partial W$. Moreover, there are indefinite metrics $(2W, g)$ and $(W, g_W)$ of signature $(2, n - 2)$ that arise from a 2-distribution with trivial Euler class and such that the restriction of $g_{W}$ to $\partial W$ yields a non-singular Lorentzian metric. 

$\bullet$ If $n \equiv 3 \mod 4$, such pseudo-Riemannian metrics exist on $2W$ and $W$ without any other assumptions. 
\end{proposition}

The proof of Proposition \ref{Proposition IDouble} is a straight-forward application of results of Reinhart \cite[Theorem 1]{[Reinhart]} and Table 1.

\section{Spin structures on pseudo-Riemannian cobordisms and their structure groups}\label{Section Structure Groups}

Let us justify now the notation $\Spin(2n,n-1)_0$ in our definition. Let $X$ be a smooth orientable compact $n$-manifold. If $X$ admits a a pseudo-Riemannian metric $(X,g)$ of signature $(p,q)$ with $p+q=n$, then it is well known that the structure group of its tangent bundle $TX$ can be reduced to $O(p,q)$ \cite{[Steenrod]}, where $O(p,q)$ is the indefinite orthogonal group of signature $(p,q)$; i.e.
\begin{align}
    O(p,q)= \{A\in \mathbb R^{n\times n}\mid \; A^TI_{p,q}A=I_{p,q}\};
\end{align}
and $I_{p,q}$ is the diagonal matrix with the first $p$ diagonal entries equal $-1$ and the last $q$ diagonal entries equal $+1$. It is easy to see that if $A\in O(p,q)$, then $\det A=\pm 1$. In particular, if $(X,g)$ is orientable, the structure group can be further reduced to $SO(p,q)$, namely to the group of indefinite orthogonal matrices of signature $(p,q)$ with positive determinant. Such group has two connected components \cite{[Gallier]}; the connected component of the identity is denoted by $SO(p,q)_0$. One can prove that the structure group of the tangent bundle $TX$ can be reduced to $SO(p,q)_0$ if $(X,g)$ is time-orientable \cite{[May]},\cite{[Gallier]}. As for the definite case, we have a double cover of $SO (p,q)_0$ by the $\Spin(p,q)_0$ group (See \cite{[Lawson]} for the definitions) which satisfies a short exact sequence of groups
\begin{align}\label{eq: short-exact-spin}
    0\to \Z_2\to \Spin(p,q)_0\xrightarrow{\text{Ad}} SO(p,q)_0 \to 0.
\end{align}

\begin{definition}
    Let $X$ be a smooth orientable compact $n$-manifold with a time-oriented pseudo-Riemannian metric $g$ of signature $(p,q)$, where $p+q=n$. Let $\pi_{SO}:F(X)\to X$ be the principal $SO(p,q)_0$-bundle of oriented orthonormal frames of $(X,g)$. A Spin-structure on $X$ is a principal $\Spin(p,q)_0$-bundle $\pi_{\Spin}:\tilde F(X)\to X$  with a 2-fold cover $\Lambda:\tilde F(X)\to F(X)$ such that the diagram
\begin{equation}\label{diagram}
     \begin{tikzcd}
        \Spin(p,q)_0\arrow[r]\arrow[d,"\text{Ad}"] & \tilde F(X)\arrow[r]\arrow[d,"\Lambda"] & X \arrow[d,"\text{Id}"]\\
        SO(p,q)_0 \arrow[r] & F(X)\arrow[r,"\pi_{SO}"] & X
    \end{tikzcd}
\end{equation}
commutes.
\end{definition}
The existence of a $\Spin(p,q)_0$-structure under these hypotheses does not depend on the metric. The obstruction for it is merely topological, as the following result exhibits.

\begin{lemma}
    Let $X$ be a compact orientable smooth $n$-manifold $X$ with a time-orientable pseudo-Riemannian metric of signature $(p,q)$, with $p+q=n$. Then $X$ admits a $\Spin(p,q)_0$-structure if and only if it admits a $\Spin(n)$-structure.
\end{lemma}
\begin{proof}
    Karaubi showed in \cite[Proposition 1.1.26]{[Kar]} that, given a splitting $TX=\xi\oplus \eta$ induced by a pseudo-Riemannian metric of signature $(p,q)$, the principal $SO(p,q)_0$-bundle $F(X)$ defined above admits a lift to $\tilde F(X)$ as in (\ref{diagram}) if and only if the following conditions are satisfied:
    \begin{align}
        w_1(\xi)+w_1(\eta)=&0 & \text{and}& & w_2(\xi)+w_2(\eta)=0.
    \end{align}
    Since $X$ is assumed to be orientable, the first equation above is immediately satisfied. Since $w_2(X)=w_2(\xi)+w_2(\eta)+w_1(\xi)\smile w_1(\eta)$ by elementary properties of Stiefel-Whitney classes and $(X,g)$ is time-orientable, we have $w_1(\xi)=w_1(\eta)=0$ and thus $w_2(X)=0$ if and only if $w_2(\xi)+w_2(\eta)=0$.
\end{proof}

% In Section \ref{Section Indefinite Metrics} we have seen that, given an orientable $n$-manifold equipped with a non-singular pseudo-Riemannian metric of signature $(2,n-2)$, its tangent bundle splits as the direct sum of two oriented subbundles of rank $2$ and $n-2$ respectively. This implies that the structure group of its frame bundle can be reduced to $\text{SO}(2,n-2)_0$, namely to the connected component of the generalized special orthogonal group $\text{SO}(2,n-2)$ containing the identity matrix.  Moreover, if we ask our manifold to be $\Spin$, we get a further reduction of the structure group to $\text{Spin}(2,n-2)_0$. For further details on such groups see \cite{[Lawson]} or \cite{[SmirnovTorres]} for the analogue in the $\Spin$-Lorentzian case.

\begin{corollary}
    Let $((W;M_1,M_2),g)$ be an $SO(2,n-1)_0$-pseudo-Riemannian cobordism. Then $((W;M_1,M_2),g)$ is a $\spcob$ if and only if $W$ admits a $\Spin(2,n-1)_0$-structure.
\end{corollary}

  \section{Spin pseudo-Riemannian cobordism groups}\label{Section Groups}

Let us discuss the definition of the $\Spin(2,n-1)_0$-pseudo-Riemannian cobordism groups that already appear in the statement of Theorem \ref{Theorem Extension2}. For each integer $n\geq 1$, let $\mathcal A_n$ be the set of diffeomorphism classes of closed $\Spin$ $n$-manifolds $M$ with the property that all of their connected components have vanishing Euler characteristic. Define following relation in $\mathcal A_n$:
    \begin{align}\label{eqn:  rel spin cob}
        M_1\sim& M_2 &\text{if and only if}& & \{M_1,M_2\}\text{ are $\Spin(2,n-1)_0$-cobordant.}
    \end{align}

\begin{proposition}\label{prop: equiv rel spin cob}
    The relation (\ref{eqn:  rel spin cob}) defines an equivalence relation on $\mathcal A_n$.
\end{proposition}

Lemma \ref{Lemma 2Distribution} is a key ingredient in the proof of Proposition \ref{prop: equiv rel spin cob}, as it gives us tools to glue different $\Spin(2,n-1)_0$-pseudo-Riemannian cobordisms with mutual boundary connected components; see \cite{[May]} for details.

\begin{definition}\label{Definition Cobordism Groups}
  The $\Spin(2,n-2)_0$-pseudo-Riemannian cobordism group is the set of the equivalence classes of the relation \ref{eqn:  rel spin cob} equipped with the disjoint union as group product and it is denoted by $\Omega^{\Spin_0}_{2,n-1}$.
\end{definition}

The reader might have already noticed that the group operation in Definition \ref{Definition Cobordism Groups} cannot be the connected sum of $M_1$ and $M_2$ as it is the case in cobordisms of other flavors. Notice that if $M_1$ and $M_2$ are even-dimensional Lorentzian manifolds, their connected sum $M_1\#M_2$ does not admit a non-singular Lorentzian metric. 

The main result of this section is the following theorem.

\begin{theorem}\label{thm: spin cob groups}
    The $\Spin(2,n-1)_0$-pseudo-Riemannian groups $\Omega^{\Spin_0}_{2,n-1}$ are abelian groups. Moreover, the Cartesian product of manifolds yields a graded ring structure
    \begin{align}\label{eqn: spin cob ring}
        \Omega^{\Spin_0}_{2,*-1}=\bigoplus_{n=1}^\infty \Omega^{\Spin_0}_{2,n-1}.
    \end{align}
\end{theorem}

The proof of Theorem \ref{thm: spin cob groups} is a straightforward adaptation of the classical arguments in cobordism theory \cite{[Milnor]},\cite{[Kirby]},\cite{[Reinhart]},\cite{[Stong]},\cite{[Yodzis]}, and a detailed proof can be found in \cite{[May]}. We refer to the ring 
 (\ref{eqn: spin cob ring}) as the $\Spin(2,n-1)_0$-pseudo-Riemannian cobordism ring.

\section{Comparison of pseudo-Riemannian cobordisms}\label{Section Comparison} In this section, we draw a comparison between the pseudo-Riemannian cobordisms of Definition \ref{Definition Kleinian Cobordism} and $\Spin(1, n - 1)_0$-Lorentzian cobordisms. The contrast between these objects sheds light on the topological restrictions imposed by the coexistence of the $\Spin$-structure with the pseudo-Riemannian structure of the cobordism. We first recall the definition of a $\Spin(1, n - 1)_0$-Lorentzian cobordism.

\begin{definition}\label{Definition Lorentzian Cobordism}
A Lorentzian cobordism between closed smooth $n$-manifolds $M_1$ and $M_2$ is a pair

\begin{align}\label{Lorentzian Cobordism}
((\hat{W}; M_1, M_2), g^L)
\end{align} 
that consists of

(A) a cobordism $(\hat{W}; M_1, M_2)$,

(B.2) a non-singular Lorentzian metric $(\hat{W}, g^L)$ with a time-like line field $V$,

(C.2) and the boundary $\partial \hat{W} = M_1 \sqcup M_2$ is space-like, i.e. $(M_1, g^L|_{M_1})$ and $(M_2, g^L|_{M_2})$ are Riemannian manifolds, where $g^L|_{M_i}$ is the restriction of $g^L$ to $M_i$.

If the cobordism $(\hat{W}; M_1, M_2)$ of Item (A) is a $\Spin$-cobordism, we say that the cobordism (\ref{Lorentzian Cobordism}) is a $\Spin(1, n - 1)_0$-Lorentzian cobordism.

 \end{definition}
 
We keep the discussion at a three-dimensional level for the sake of brevity, although similar comparisons apply to any dimension. More precisely, we pivot the comparison and the discussion of the topological restrictions on the following result.

\begin{theorem}\label{Theorem Lorentzian}Gibbons-Hawking \cite{[GibbonsHawking1]}, Smirnov-Torres \cite{[SmirnovTorres]}. Let $\{M_1,  M_2\}$ be closed oriented 3-manifolds. The following conditions are equivalent
\begin{enumerate}
\item There exists a $\Spin(1, 3)_0$-Lorentzian cobordism\begin{center}$((\hat{W}; M_1, M_2), g),$\end{center}where $(\hat{W}; M_1, M_2)$ is a $\Spin$-cobordism 
\item $\hat{\chi}_{\Z/2}(M_1) = \hat{\chi}_{\Z/2}(M_2)$
\item $\hat{W}$ is parallelizable and $\chi(\hat{W}) = 0$.
\end{enumerate}

There is a group isomorphism\begin{align}\label{4D Isomorphism Lorentzian}\Omega^{\Spin_0}_{1, 3} \rightarrow \Z/2.\end{align}
\end{theorem}

While the existence of a $\Spin$-cobordism imposes no restrictions on the boundary 3-manifold since the third cobordism group $\Omega^\Spin_3$ is trivial \cite{[Milnor]}, the presence of the required Lorentzian metric on $\hat{W}$ forces its Euler characteristic to be $\chi(\hat{W}) = 0$ \cite{[Reinhart]}. For the latter structure to coexist with the $\Spin$-cobordism, the Kervaire semi-characteristic of the boundary 3-manifolds must coincide as indicated by Theorem \ref{Theorem Relation}. This invariant gives us the isomorphism (\ref{4D Isomorphism Lorentzian}). 

The corresponding statement for $\Spin(2, 2)_0$-pseudo-Riemannian cobordisms is the following.

\begin{theorem}\label{Theorem Kleinian} Let $\{M_1,  M_2\}$ be closed oriented 3-manifolds. The following conditions are equivalent
\begin{enumerate}
\item There exists a $\Spin(2, 2)_0$-pseudo-Riemannian cobordism \begin{center}$((W; M_1, M_2), g),$\end{center}
\item $\hat{\chi}_{\Z/2}(M_1) = \hat{\chi}_{\Z/2}(M_2)$
\item $W$ is parallelizable.
\end{enumerate}

There is a group isomorphism\begin{align}\label{4D Isomorphism}\Omega^{\Spin_0}_{2, 2} \rightarrow \Z/2.\end{align}
\end{theorem}

The proof of Theorem \ref{Theorem Kleinian} is given in Section \ref{Section Proof of Theorem Kleinian}. The reader will notice that while both 4-manifolds $\hat{W}$ and $W$ in Theorem \ref{Theorem Lorentzian} and Theorem \ref{Theorem Kleinian} are parallelizable, the Euler characteristic of the cobordism $W$ need not be zero. The following example is illustrative of the situation.

\begin{example}\label{Example Both Cobordisms} The product of a 2-disk with the round 2-sphere admits a $\Spin$-structure as well as an indefinite metric of signature (2, 2)\begin{align}\label{Cob Example}(D^2, -dr^2 - r^2d \theta^2)\times (S^2, g_{S^2})\end{align}that restrict to a $\Spin$-structure and a Lorentzian metric on the boundary\begin{align}\label{Boundary Example}(S^1, - d\theta^2) \times (S^2, g_{S^2})\end{align}as indicated in Example \ref{Example A}. The Euler characteristic of (\ref{Cob Example}) is $\chi(D^2\times S^2) = 2$ and it does not admit a Lorentzian metric that restricts to a Riemannian metric on (\ref{Boundary Example}). The connected sum $$\hat{W} = D^2\times S^2\#S^1\times S^3,$$ on the other hand, does support both kinds of non-singular pseudo-Riemannian metrics as well as $\Spin$-structures that restrict to (\ref{Boundary Example}).
\end{example}

The phenomenon displayed in Example \ref{Example Both Cobordisms} occurs for all oriented 3-manifolds. 

\begin{corollary}Let $M_1$ and $M_2$ be closed oriented 3-manifolds. There is a $\Spin(1, 3)_0$-Lorentzian cobordism if and only if there is a $\Spin(2, 2)_0$-pseudo-Riemannian cobordism.
\end{corollary}

We now elucidate on the reason behind the difference in the obstructions. The existence of a 2-distribution is not equivalent to the existence of a pair of linearly independent and nowhere vanishing vector fields as shown in the work of Atiyah \cite{[Atiyah]}, Frank \cite{[Frank]}, Hirzebruch-Hopf \cite{[HirzebruchHopf]}, Matsushita \cite{[Matsushita2]}, Thomas \cite{[Thomas1], [Thomas2], [Thomas6], [Thomas7], [Thomas3], [Thomas4], [Thomas5]}. The existence of a pair of linearly independent vector fields on a manifold requires for its Euler characteristic to vanish in the even-dimensional case and for its Kervaire semi-characteristic to vanish in the odd-dimensional case.

As we end this section, we take the opportunity to amend \cite[Corollary E]{[SmirnovTorres]}. The correct statement is as follows. 

\begin{corollary}Let $M_1$ and $M_2$ be closed smooth Spin 4-manifolds. There is a $\Spin(1, 4)_0$-Lorentzian cobordism $((W; M_1, M_2), g^L)$ if and only if $\chi(M_1) = \chi(M_2)$ and $\sigma(M_1) = \sigma(M_2)$.
\end{corollary}

\section{Proofs}\label{Section Proofs}

The proofs of the results that are mentioned in the introduction have the following structure.

\subsection{Proof of Theorem \ref{Theorem Kleinian}}\label{Section Proof of Theorem Kleinian}We first show the equivalence (1) $\Leftrightarrow$ (2) and begin with (1) $\Leftarrow$ (2). Suppose $M_1$ and $M_2$ are two closed oriented 3-manifold whose Kervaire semi-characteristics satisfy $\hat{\chi}_{\Z/2}(M_1) = \hat{\chi}_{\Z/2}(M_1)$. As it was mentioned in the previous section, the third $\Spin$-cobordism group is $\Omega^\Spin_3 = \{0\}$ and there is a $\Spin$-cobordism $(\hat W; M_1, M_2)$.  Theorem \ref{Theorem Relation} implies that\begin{align}\label{4D Implication}\chi(\hat W) + \hat{\chi}_{\Z/2}(\partial \hat W) \equiv 0 \mod 2.\end{align}Thus, we have that\begin{align}\chi(\hat W) \equiv 0\mod 2.\end{align} Take connected sums of $\hat{W}$ with copies of $S^1\times S^3$ and $S^2\times S^2$ and obtain a manifold  $$W:=\hat W\consum k_1(S^1 \times S^3)\consum k_2(S^2\times S^2)$$ which has zero Euler characteristic, by choosing $k_1$ and $k_2$ appropriately. Proposition \ref{Proposition IDouble} allows us to conclude the proof of the implication. To show that the implication (1) $\Rightarrow$ (2) holds, we proceed as follows. Assume that there is a $\Spin(2, 2)_0$-pseudo-Riemannian cobordism $((W; M_1, M_2), g)$ and consider the closed smooth 4-manifold with the indefinite metric $(2W, g_{2W})$ given as the double of $(W, g)$. Applying Theorem \ref{Theorem Main} to $2W$ we get that $\chi(W) = 0 \mod 2$ and the Kervaire semi-characteristics of $M_1$ and $M_2$ coincide by (\ref{4D Implication}). We conclude that the equivalence (1) $\Leftrightarrow$ (2) holds. The implication (1) $\Rightarrow$ (3) follows from $W$ being an almost parallelizable manifold with non-empty boundary. Such a manifold is stably-parallelizable, and hence parallelizable  \cite{[KervaireMilnor]}, \cite[\S 7, \S 8]{[Kosinski]}. The implication (3) $\Rightarrow$ (2) follows from Theorem \ref{Theorem Relation}. 

The Kervaire semi-characteristic yields an isomorphism\begin{align}\hat{\chi}_{\Z/2}:\Omega^{\Spin}_{2, 2}\rightarrow \Z/2\end{align} and its generator is $S^3$; see \cite[Theorem C]{[SmirnovTorres]}.

\hfill $\square$

\subsection{Proof of Theorem \ref{Theorem Extension3}}
Let's begin by studying the case $n\equiv 0 \mod 2$ and let $(W;M_1,M_2)$ be a Spin-cobordism with $\chi(M_1)=0=\chi(M_2)$. Since 
$$2\chi(W)=\chi(2W)+\chi(\partial W)=0$$ 
by Lemma \ref{Lemma Characteristic Double}, Proposition \ref{Proposition IDouble} implies that there is a metric $g$ on $W$ for which $((W;M_1,M_2),g)$ is a $\Spin(2, n - 1)_0$-pseudo-Riemannian cobordism. Conversely, if $((W;M_1,M_2),g)$
is a $\spcob$, then the fact that $g$ restricts to a Lorentzian metric on $\partial W=M_1\sqcup M_2$ implies that $\chi(M_1)=\chi(M_2)=0.$

We now address the cases $n\equiv 1, 3, 5 \mod 8$.  If $\hat{\chi}_{\Z/2}(M_1) = \hat{\chi}_{\Z/2}(M_2) \mod 2$, there is a $\Spin(1,n)_0$-Lorentzian cobordism $((W;M_1, M_2),g^L)$ by  \cite[Theorem D]{[SmirnovTorres]}. In particular, we have that $$\chi(2W)=2\chi(W)=0.$$ Being $M_1, M_2$ odd-dimensional, we also have that $$\chi(M_1)=0=\chi(M_2)$$ and hence we can apply Proposition \ref{Proposition IDouble} to obtain a $\Spin(2,n-1)_0$-pseudo-Riemannian cobordism $((W;M_1, M_2),g)$. Suppose now that $((W;M_1, M_2),g)$ is a $\Spin(2,n-1)_0$-pseudo-Riemannian cobordism and let us restrict to the case $n\equiv 3\mod 8$. We have that $\chi(W)\equiv 0 \mod 2$ by Theorem \ref{Theorem Main}. Theorem \ref{Theorem Relation} allows us to conclude that $\hat{\chi}_{\Z/2}(M_1) = \hat{\chi}_{\Z/2}(M_2) \mod 2$.

Consider the case $n\equiv 7 \mod 8$. By \cite[Theorem D]{[SmirnovTorres]} there is a $\Spin(1, n)_0$-Lorentzian cobordism $((W ; M_1, M_2), g^L)$. A result of Reinhart \cite{[Reinhart]} says that $\chi(W)=0$. Hence, Proposition \ref{Proposition IDouble} implies the existence of a $\Spin(2,n-1)_0$-pseudo-Riemannian cobordism $((W;M_1,M_2),g)$.

\hfill $\square$

\begin{remark}\label{Remark Complete} The reason for the abscence of a complete characterization in the case $n \equiv 1 \mod 4$ of Theorem \ref{Theorem Extension3} is essentially due to the lack of existence results in literature of distributions of tangent $2$-planes on closed $4q+2$-dimensional manifolds. In particular, it is not clear whether the existence of a $\spcob$ implies the existence of a $\Spin(1,n)_0$-Lorentzian cobordism. Thus, the missing implication.
    By inspecting the proof of Theorem \ref{Theorem Extension3}, one realizes that solving this issue is equivalent to giving an answer to the following open question.
    \end{remark}
    \begin{center}
        \textbf{Question:} Is it possible to find two closed $n$-manifolds $M_1, M_2$ which are $\Spin(2,n-1)_0$-pseudo-Riemannian cobordant and such that the Euler class of the $2$-distribution inducing the pseudo-Riemannian metric on the cobordism is necessarily non-trivial?
    \end{center}

On the other hand, from the proof of Theorem \ref{Theorem Extension3} we get that we may ask for $\text{Span }W\geq 2$ for a $\spcob$ whenever $n \not \equiv 1 \mod 4$. This yields the following corollary.

\begin{cor}
    Let $n$ be positive integer with $n \not\equiv 1 \mod 4$. Let $\{M_1,M_2\}$ be two smooth closed $\Spin (2,n-1)_0$-pseudo-Riemannian cobordant $n$-manifolds. There exists a smooth Spin $(n+1)$-manifold $W$ and two everywhere linearly independent vector fields $X,Y\in \mathfrak X(W)$ such that 
    \begin{itemize}
        \item $(W;M_1,M_2)$ is a $\Spin$-cobordism;
        \item $X$ is interior normal to $M_1$ and exterior normal to $M_2$.
    \end{itemize}
\end{cor}

In particular, the existence of the $\Spin$-cobordism $(W;M_1,M_2)$ and a nowhere vanishing vector field  $X\in \mathfrak X(W)$ as above is equivalent to the existence of a $\Spin(1,n)_0$-Lorentzian cobordism $((W;M_1,M_2),g)$ (see \cite{[SmirnovTorres]}). 

\subsection{Proof of Theorem \ref{Theorem Extension2}}\label{Section Proof of Theorem Extension2} The three-dimensional case has been addressed in Theorem \ref{Theorem Kleinian}. We argue the four-dimensional case first; cf. \cite[Proof of Corollary E]{[SmirnovTorres]}. Suppose $((W; M_1, M_2), g)$ is a $\Spin(2, 3)_0$-pseudo-Riemannian cobordism. Since $(W; M_1, M_2)$ is a $\Spin$-cobordism, it is in particular an oriented cobordism and therefore we have the condition $\sigma(M_1)=\sigma(M_2)$. Moreover, the existence of a Lorentzian metric on $M_i$ for $i = 1, 2$ implies that $\chi(M_i) = 0$. To prove the converse, we argue as follows. If $M_1$ and $M_2$ are closed $\Spin$ 4-manifolds with the same signature, there is a $\Spin$-cobordism $(W; M_1, M_2)$ \cite{[Thom]}, \cite[Chapter VIII]{[Kirby]}. The existence of an indefinite $(2, 3)$-metric on $W$ restricting to a Lorentzian one on the boundary follows from the vanishing of the Euler characteristics of $M_1$ and $M_2$ and Theorem \ref{Theorem Extension3}. The group isomorphism $$\Omega^{\Spin}_{2,3}\to \Z$$ is given by the signature. In particular, the map $[M]\mapsto \sigma(M)$ is well defined, being the signature an oriented cobordism invariant.

Let us address now the six-dimensional case. The existence of a $\Spin(2, 5)_0$-pseudo-Riemannian cobordism $((W; M_1, M_2), g)$ implies that the Euler characteristics of $M_1$ and $M_2$ are $\chi(M_i) = 0$ for $i = 1, 2$ \cite{[Thom]}. The converse follows from results of Milnor and Thomas. Since the sixth cobordism group is $\Omega^{\Spin}_6 = \{0\}$ \cite{[Milnor]}, we know that there is a $\Spin$-cobordism $(W; M_1, M_2)$. The conclusion follows again from Theorem \ref{Theorem Extension3}. The group $\Omega^{\Spin_0}_{2, 5}$ is hence trivial.

In the seven-dimensional case, Milnor observed that $\Omega^{\Spin}_7 = \{0\}$ \cite[p. 201]{[Milnor]} and any two closed $\Spin$ 7-manifolds bound a $\Spin$ 8-manifold. The existence of a $\Spin(2, 6)_0$-pseudo-Riemannian cobordism follows from Theorem \ref{Theorem Extension3}, and the group $\Omega^{\Spin_0}_{2, 6}$ is hence trivial.

\hfill $\square$


\begin{thebibliography}{99}
\bibitem{[AltyChamblin]} L. J. Alty and A. Chamblin, \emph{Spin structures on Kleinian manifolds}, Class. Quantum Grav. 11 (1994), 2411 - 2415.

\bibitem{[Atiyah]} M. F. Atiyah, \emph{Vector fields on manifolds}, Arbeitsgemeinschaft f\"ur Forschung des Landes Nordhein-Westfalen, Heft 200

%\bibitem{[BredonKosinski]} G. E. Bredon and A. A. Kosinski, \emph{Vector fields on $\pi$-manifolds}, Ann. of Math. 84 (1966), 85 - 90.

\bibitem{[Chamblin]} A. Chamblin, \emph{Some applications of differential topology in general relativity}, J. Geom. Phys. 13 (1994), 357 - 377.

%\bibitem{[DoldWhitney]} A. Dold and H. Whitney, \emph{Classification of oriented sphere bundles over a 4-complex}, Ann. of Math. 69 (1959), 667 - 677.

\bibitem{[Frank]} D. Frank, \emph{On the index of a tangent 2-field}, Topology 11 (1972), 245 - 252.


\bibitem{[Gallier]} J. Gallier, \emph{Cliffor Algebras, Clifford Groups, and a Generalization of the Quaternions} (2014), arXiv:0805.0311v3


\bibitem{[Geiges]} H. Geiges, \emph{An introduction to contact topology}, Cambridge Studies in Advanced Mathematics, vol. 109, Cambridge University Press, Cambridge, 2008. 

\bibitem{[Geroch]} R.P. Geroch, \emph{Topology in general relativity}, J. Mathematical Phys. 8 (1967), 782 - 786.

\bibitem{[GibbonsHawking1]} G. W. Gibbons and S. W. Hawking, \emph{Selection rules for topology change}, Commun. Math. Phys. 148 (1992), 345 - 352.

\bibitem{[Hirsch]} M. W. Hirsch, \emph{Differential Topology}, Springer-Verlag New York, Berlin Heidelberg London Paris Tokyo Hong Kong Barcelona Budapest, 1994

\bibitem{[HirzebruchHopf]} F. Hirzebruch and H. Hopf, \emph{Felder von Fl\"achenelementen in 4-dimensionalen Mannigfaltigkeiten}, Math. Ann. 136 (1958), 156 - 172.

%\bibitem{[Ibarra]} D. Ibarra, \emph{Vector fields on orientable 7-manifolds}, arXiv:2207.12149v2, 2022.

\bibitem{[JohnsonWalton]} F. E. A. Johnson and J. P. Walton, \emph{Parallelizable manifolds and the fundamental group}, Mathematika 47 (2000), 165 - 172. 

\bibitem{[Kar]} M. Karoubi,  \textit{Algèbres de Clifford et $K$-théorie}. Annales scientifiques de l'École Normale Supérieure, Serie 4, Volume 1 (1968) no. 2, pp. 161-270. doi : 10.24033/asens.1163

\bibitem{[Kervaire]} M. Kervaire, \emph{Courbure int\'egrale g\'en\'eralis\'ee et homotopie} (French), Math. Ann. 131 (1956), 219 - 252.

\bibitem{[KervaireMilnor]} M. Kervaire and J. Milnor, \emph{Groups of homotopy spheres: I}, Ann. of Math. 77 (1963), 504 - 537.

\bibitem{[Kirby]} R. C. Kirby, \emph{The topology of 4-manifolds}, Lect. Notes in Math. 1374, Springer-Verlag, Berlin, 1989.

\bibitem{[Kosinski]} A. A. Kosinski, \emph{Differential manifolds}, Academic Press, Inc. 1993.

\bibitem{[Lawson]} H. B. Lawson and M. Michelsohn, \emph{Spin Geometry}, Princeton University Press, 1989.

\bibitem{[LusztigMilnorPeterson]} G. Lusztig, J. Milnor and F. P. Peterson, \emph{Semi-characteristics and cobordism}, Topology 8 (1969), 357 - 359.

%\bibitem{[MantioneTorres]} A. Mantione and R. Torres, \emph{Geography of 4-manifolds with positive scalar curvature}, Exp. Math. 39 (2021), 566 - 582. 

\bibitem{[Matsushita1]} Y. Matsushita, \emph{Fields of 2-planes on compact simply-connected smooth 4-manifolds}, Math. Ann. 280 (1988), 687 - 689.

\bibitem{[Matsushita2]} Y. Matsushita, \emph{Fields of 2-planes and two kinds of almost complex structures on compact 4-dimensional manifolds}, Math. Z. 207 (1991), 281 - 291.

\bibitem{[May]} V. G. May Custodio, \emph{Unpublished Master's Thesis: On the $\Spin(2,n-1)_0-$pseudo-Riemannian cobordism groups}. Universit\'a degli Studi di Trieste/SISSA (2023).

\bibitem{[Milnor3]} J. Milnor, \emph{Remarks concerning spin manifolds}, Differential and Combinatorial Topology: A Symposium in Honor of Marston Morse, Princeton: Princeton University Press, (1965), pp. 55-62. 

\bibitem{[Milnor]} J. Milnor, \emph{Spin structures on manifolds}, L'Enseignement Math\'ematique, 9 (1963), 198 - 203.

\bibitem{[Milnor2]} J. Milnor, \emph{Topology from the differentiable viewpoint}, The University Press of Virginia, Charlottesville




\bibitem{[ONeill]} B. O'Neill, \emph{Semi-Riemannian geometry with applications to relativity}, Academic Press, Inc. 1983.

\bibitem{[Reinhart]} B. L. Reinhart, \emph{Cobordism and the Euler number}, Topology 2 (1963), 173 - 177.

\bibitem{[Scorpan]} A. Scorpan, \textit{The Wild World of 4-Manifolds}, American Mathematical Soc. (2005), ISBN: 0821837494

\bibitem{[SmirnovTorres]} G. Smirnov and R. Torres, \emph{Topology change and selection rules for high-dimensional $\Spin(1, n)_0$-Lorenzian cobordisms}, Trans. Amer. Math. Soc. 373 (2020), 1731 - 1747.

\bibitem{[Sorkin]} R. D. Sorkin, \emph{Topology change and monopole creation}, Phys. Rev. Lett. D 33 (1986), 978 - 982.

\bibitem{[Steenrod]} N. Steenrod, \emph{The topology of fibre bundles}, Princeton Math. Series, vol. 14, Princeton University Press, Princeton, N. J., 1951.

\bibitem{[Stong]} R. E. Stong, \emph{Notes on cobordism theory}, Mathematical notes, Princeton University Press, Princeton, N.J., University of Tokyo Press, Tokyo, 1968.

\bibitem{[Thom]} R. Thom, \emph{Quelque propri\'et\'es globales des vari\'et\'es diff\'erentiables}, Comm. Math. Helv. 28 (1954), 17 - 86.

\bibitem{[Thomas6]} E. Thomas, \emph{Seminar on Fiber Spaces}, Lect. Notes in Math. 13, Springer-Verlag, Berlin, 1966.

\bibitem{[Thomas1]} E. Thomas, \emph{The index of a tangent 2-field}, Comment. Math. Helv. 42 (1967), 86 - 110.

\bibitem{[Thomas2]} E. Thomas, \emph{Fields of tangent 2-planes on even-dimensional manifolds}, Ann. of Math. 86 (1967), 349 - 361. 

\bibitem{[Thomas3]} E. Thomas, \emph{Fields of tangent k-planes on manifolds}, Invent. Math. 3 (1967), 334 - 347. 

\bibitem{[Thomas7]} E. Thomas, \emph{Postnikov invariants and higher order cohomology operations}, Ann. of Math. 85 (1967), 184 - 217.

\bibitem{[Thomas4]} E. Thomas, \emph{Vector fields on low dimensional manifolds}, Math. Z. 103 (1968), 85 - 93.

\bibitem{[Thomas5]} E. Thomas, \emph{Vector fields on manifolds}, Bull. Amer. Math. Soc. 75 (1969), 643 - 683.

\bibitem{[Thomas8]} E. Thomas, \emph{Some remarks on vector fields}, H-spaces, Actes de la r\'eunion de Neuchatel (Suisse), Lecture Notes in Math. 196 (1970), 107 - 113. 

\bibitem{[Thomas9]} E. Thomas, \emph{Real and complex vector fields on manifolds}, J. Math, and Meeh. 16
(1967) 1183-1206.

\bibitem{[Yodzis]} P. Yodzis, \emph{Lorentz cobordism}, Comm. Math. Phys. 26 (1072), 39 - 52. 

\bibitem{[MEZadeh]} M. E. Zadeh, \emph{On cut-and-paste invariant of Kervaire semi-characteristic}, arXiv:1110.2447v1, 2011.


\end{thebibliography}
\end{document}